\documentclass[11pt]{amsart}

\usepackage{geometry}
\usepackage{amsmath, amssymb, amsfonts, amsopn, amsthm}

\usepackage{makecell, longtable, booktabs, pifont} 
\usepackage[parfill]{parskip} 
\usepackage{graphicx, tikz} 
\usepackage{url, float} 
\usepackage{subfigure} 
\usepackage[percent]{overpic} 
\usepackage{extarrows}

\definecolor{blu}{gray}{0.7}  

\newcommand{\knot}[1]{\includegraphics[page=#1,scale=0.8]{knot}}

\newcommand{\N}{\mathbb{N}} \newcommand{\Z}{\mathbb{Z}} \newcommand{\R}{\mathbb{R}} 

  \newcommand{\rpt}{\mathbb{R}\text{P}^3}
\newcommand{\lpq}{L(p,q)} \newcommand{\lpo}{L(p,1)}

\newcommand{\calS}{\mathcal{S}}\newcommand{\calB}{\mathcal{B}}
\newcommand{\calLfr}{\mathcal{L_\text{fr}}}
\newcommand{\calLor}{\mathcal{L_\text{or}}}

\newcommand{\df}[1]{\emph{#1}}
\newcommand{\gauss}[1]{\begin{center}{\tt #1}\end{center}}
\newcommand{\rob}{\partial}

\newcommand{\ri}{\Omega_1} \newcommand{\rii}{\Omega_2} \newcommand{\riii}{\Omega_3}    
\newcommand{\slajd}{\textit{SL}_{p,q}}
\newcommand{\gc}{\mathcal{W}}

\newcommand{\shsm}{\mathcal{S}_3}
\newcommand{\skbsm}{\mathcal{S}_{2,\infty}}

\newcommand{\kbsm}{\mathit{KBSM}}
\newcommand{\hsm}{\mathit{HSM}}
\newcommand{\GCD}{\mathit{GCD}}

\newcommand{\wind}{\mathit{wind}}

\newtheorem{thm}{Theorem}[section]

\newtheorem{cor}[thm]{Corollary}

\newtheorem{prop}[thm]{Proposition}
\theoremstyle{definition}

\theoremstyle{remark}
\newtheorem{remark}[thm]{Remark}
\newtheorem*{example}{Example}
\numberwithin{equation}{section}

\begin{document}



	
	
	

	\title[Tabulation of Prime Knots in Lens Spaces]
	{Tabulation of Prime Knots in Lens Spaces}
	
	\author[B. Gabrov\v sek]{Bo\v stjan Gabrov\v sek}
	
	\address{%
		University of Ljubljana\\
		Faculty of Mathematics and Physics\\
		Jadranska ulica 19\\
		1000 Ljubljana\\Slovenia}
	
	\email{bostjan.gabrovsek@fmf.uni-lj.si}
	
	\subjclass{Primary 57M27; Secondary 57M50}
	
	\keywords{knot table, tabulation, classification, solid torus, lens space, skein module, hyperbolic structure}
	
	\date{\today}

\begin{abstract}
\noindent Using computational techniques we tabulate prime knots up to five crossings in the solid torus and the infinite family of lens spaces $\lpq$. For these knots we calculate the second and third skein module and establish which prime knots in the solid torus are amphichiral. Most knots are distinguished by the skein modules. For the handful of cases where the skein modules fail to detect inequivalent knots, we calculate and compare the hyperbolic structures of the knot complements. We were unable to resolve a handful of 5-crossing cases for $p\geq 13$.


\let\thefootnote\relax\footnote{\textit{Mathematics Subject Classification 2010:} 57M27, 57M50}
\let\thefootnote\relax\footnote{\textit{Keywords:} Knot, tabulation, skein module; solid torus, lens space.}

\end{abstract}
\maketitle
\section{Introduction}\label{sec:intro}

The first knot table (a list of possible embeddings of a simply closed curve into $S^3$ or $\R^3$ up to ambient isotopy) was published by P. G. Tait in 1884~\cite{T} and it is believed that this table represents the very beginning of the study of mathematical knots.
Knot tables grew in size throughout the centuries and are now known up to 16 crossings~\cite{HTW}.
It seems natural that knot tables should be expanded at least to the most simple of closed 3-manifolds -- lens spaces (the infinite family of closed 3-manifolds with Heegaard genus 1).
Up to now, the only such attempt has been made for the projective space $\rpt  \approx L(2,1)$ in \cite{Dr2}. In \cite{BM} a tabulation for the solid torus can be found, but the tabulation is not complete, since knots are considered only up to a so-called flip (a certain symmetry of the solid torus). Flipped knots are hard to distinguish since most knot invariants fail to detect them. In this paper we append a genuine knot table for the solid torus and the lens spaces $\lpq$, $0 < q < p$, $\GCD(p,q)=1$ to this modest list. The knots are tabulated for up to 5 crossings, which might not seem much, but the difficulty lies in the fact that we present knot tables for an infinite class of 3-manifolds.
The tabulation is complete for up to 4 crossings, there are two unresolved 5-crossing cases in $L(p,1)$ for $p\geq13$ and 5 unresolved 5-crossing cases in $L(p,q)$ for $p\geq13$ and $q\geq2$.

By the standard inclusion $i_{p,q}: T \hookrightarrow \lpq$ of the solid torus to the lens space (see Section~\ref{sec:lpq}), each knot $K: S^1 \hookrightarrow T$ in the solid torus defines a knot $i_{p,q} \circ K: S^1 \hookrightarrow \lpq$ in the lens space. Thus, for each $\lpq$ a subset of prime knots in $T$ yields the knot table in the lens space.
The method of providing the knot tables is straightforward: for the given space we generate all knot diagrams up to $n$ crossings and classify them by ambient isotopy. The minimal diagram of each class is an entry in the knot table.
The classification itself has been made by computer with detailed results and complete source code available in~\cite{code}, the algorithm is presented in Section~\ref{sec:alg} and the final results are summarized in Section~\ref{sec:results}.

\section{Knots in \boldmath{$\lpq$}}\label{sec:lpq}

We start off by defining knot diagrams used in the knot tables and overview knot invariants used to establish knot inequivalences. 

\subsection{Knot diagrams and skein modules}\label{sec:kdasm}

\begin{figure}[htb]
\centering
\subfigure[]{\begin{overpic}[page=1]{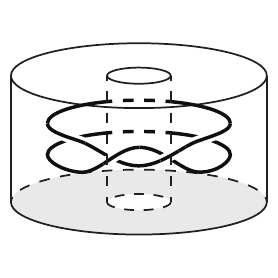}\end{overpic}}\hspace{1.0cm}\label{fig:pdda}
\subfigure[]{\begin{overpic}[page=2]{images}\end{overpic}}\hspace{1.0cm}\label{fig:pddb}
\subfigure[]{\begin{overpic}[page=3]{images}\end{overpic}}\label{fig:pddc}
\caption{Construction of a punctured disk diagram of a link in the solid torus.}
\label{fig:pdd}
\end{figure}

Let $K$ be a tame knot in the solid torus $T = A \times I$, with $A$ being an annulus (Figure~\ref{fig:pdd}(a)). By standard abuse of notation $K$ can either represent the map $S^1 \hookrightarrow T$ or the image of $S^1$ in $T$.

A \df{punctured disk diagram} of a knot $K$ is the regular projection of $K$ on $A$, keeping the information of over- and undercrossings (Figure~\ref{fig:pdd}(b)).

We resolve the inconvenience of drawing the annulus by making a dot (a puncture) in the region of $\R^2 \supset A$ that bounds the inner component of $\rob A$ and assume that the outer component of $\rob A$ lies in the unbounded region of $\R^2$ (Figure~\ref{fig:pdd}(c)). We call the dotted region the \df{0-region} and the unbounded region the \df{$\infty$-region}.

\begin{figure}[htb]
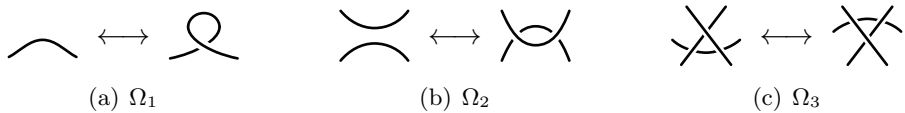

\centering
\subfigure[$\ri$]{\begin{overpic}[page=18]{images}\end{overpic}\raisebox{10pt}{$\;\;\longleftrightarrow\;\;$}\begin{overpic}[page=19]{images}\end{overpic}\label{fig:reid1}}
\hspace{1cm}
\subfigure[$\rii$]{\begin{overpic}[page=16]{images}\end{overpic}\raisebox{10pt}{$\;\;\longleftrightarrow\;\;$}\begin{overpic}[page=17]{images}\end{overpic}}
\hspace{1cm}
\subfigure[$\riii$]{\begin{overpic}[page=14]{images}\end{overpic}\raisebox{10pt}{$\;\;\longleftrightarrow\;\;$}\begin{overpic}[page=15]{images}\end{overpic}}
\caption{Classical Reidemeister moves.}
\label{fig:reid}
\end{figure}

The Reidemeister moves of a punctured disk diagram correspond to the classical Reidemeister moves $\ri$, $\rii$, and $\riii$ (Figure~\ref{fig:reid}), except that we cannot perform any move through the puncture.

\begin{figure}[htb]
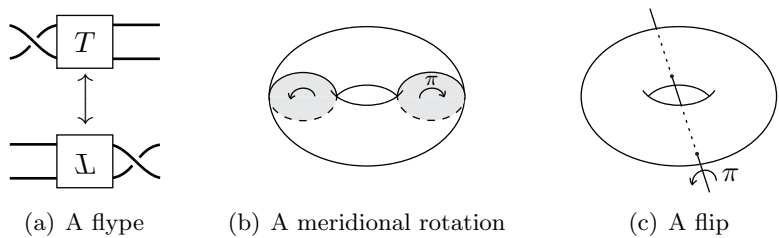

	\subfigure[A flype]{\begin{overpic}[page=38]{images}\put(36,76){$T$}\put(36,22){\scalebox{1}[-1]{$T$}}\put(40,37){\rotatebox{90}{$\longleftrightarrow$}}\end{overpic}\label{fig:flype}}
		\subfigure[A meridional rotation]{\hspace{3em}\begin{overpic}[page=44]{images}\put(77,54){$_\pi$}\end{overpic}\label{fig:mrot}\hspace{3em}}
	\subfigure[A flip]{\begin{overpic}[page=40]{images}\put(70,7){$\pi$}\end{overpic}\label{fig:flip}}
	\caption{A flype, flip, and meridional rotation .}
	\label{fig:isot}
\end{figure}

A \df{flype} is an isotopy move that consists of twisting a part of a knot using a rotation by $\pi$ as indicated in Figure~\ref{fig:flype}. A \df{meridional rotation} \cite{BM} is the self-homeomorphism of $T$ that rotates each meridional disk of $T$ by $\pi$, see Figure~\ref{fig:mrot}. If we reflect a knot $K \in T$ through $A\times\{\frac{1}{2}\}$ of $T=A\times I$, we denote the resulting knot by $\overline{K}$ and call it the \df{mirror} of $K$. 
A \df{flip} \cite{BM} is the rotation of $T$ around an axis indicated on Figure~\ref{fig:flip}. 
Note that flypes and meridional rotations are ambient isotopies in $T$, whereas a flip is not an ambient isotopy, but an orientation preserving homeomorphism of $T$.

\begin{figure}[htb]
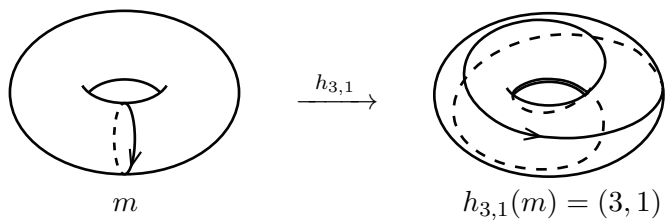

\centering
\begin{overpic}[page=4]{images}\put(45,-13){$m$}\end{overpic}
\raisebox{27pt}{$\;\;\;\;\;\xrightarrow{\;\;h_{3,1}\;\;}\;\;\;\;\;$}
\begin{overpic}[page=5]{images}\put(14,-13){$h_{3,1}(m)=(3,1)$}\end{overpic}\vspace{0.3cm}
\caption{The boundary homeomorphism $h_{3,1}:\rob T_1 \rightarrow \rob T_2$.}
\label{fig:lpq-tori}
\end{figure}

The lens space $\lpq$, $0 < q < p$, $\GCD(p,q)=1$, is the the result of glueing two solid tori $T_1$ and $T_2$ together by their boundary via the homeomorphism $h_{p,q} : \rob T_1 \rightarrow \rob T_2$ that takes the meridian of $\partial T_1$ to the $(p,q)$-curve on $\partial T_2$ (Figure~\ref{fig:lpq-tori}).

To obtain a diagram of a knot $K$ in $\lpq$ we first isotope $K$ into the first component $T_1$ and project it to the annulus $A$ of $T_1 = A\times I$. Such a diagram corresponds to the punctured disk diagram of a knot in $T_1$ (c.f.~\cite{CMM}).

\begin{figure}[htb]
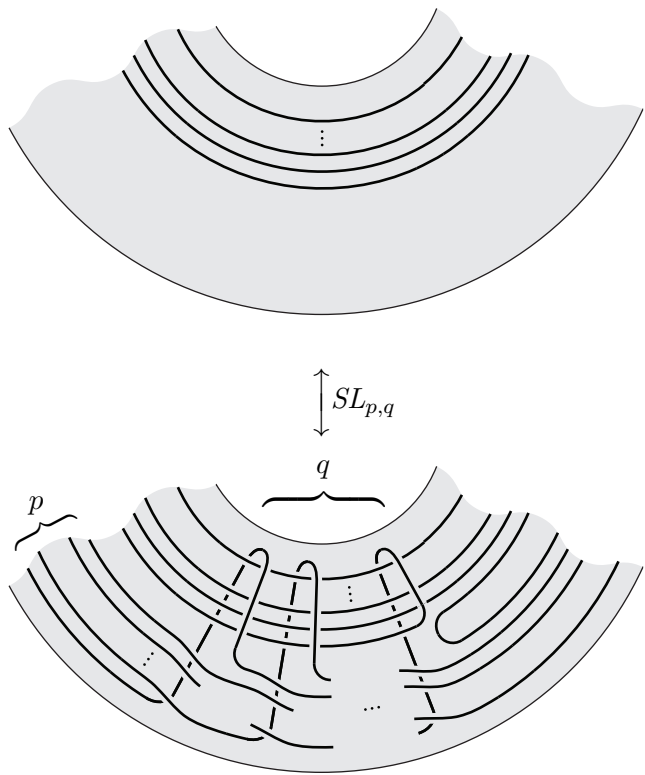

\centering
\begin{overpic}[page=23]{images}\end{overpic} \\
\hspace{0.8cm}\rotatebox{-90}{ $\xlongleftrightarrow{\rotatebox{90}{$\slajd$}}$}\\[0.3cm]
\begin{overpic}[page=22]{images}
\put(1,34){\rotatebox{118}{$\left.\begin{tabular}{c}\\[0.5cm]\end{tabular}\right\}$}}
\put(4.5,42){$p$}
\put(40,37){\rotatebox{90}{$\left.\begin{tabular}{c}\\[1.2cm]\end{tabular}\right\}$}}
\put(48.5,47.5){$q$}
\end{overpic}
\caption{The slide move for $\lpq$.}
\label{fig:slajdmove}
\end{figure}

We equip these diagrams with an additional isotopy move $\slajd$ also known as the \df{slide move}~\cite{HP1, BM4} or in some literature as the \df{band move}~\cite{LR1, LR2, La2, DL1}. This move arises from the gluing map $h_{p,q}$ and is presented in Figure~\ref{fig:slajdmove}. One can visualize the move by sliding an arc of the knot over the meridional disk of $T_1$ glued to $T_2$.


\begin{prop}[Hoste, Przytycki~\cite{HP1}]Two punctured disk diagrams represent the same link in $\lpq$ if and only if one can be transformed into the other by a finite sequence of Reidemeister moves $\ri$, $\rii$, $\riii$, and $\slajd$. \end{prop}

While there is a large amount of literarture available on the ability of knot polynomials to detect different knots in $S^3$, very little is known about $L(p,q)$-specialized polynomial invariants.
This motivates us to use skein modules as our primary invariant to distinguish inequivalent knots. Namely, we use the Kauffman bracket skein module ($\kbsm$) and the HOMFLYPT skein module ($\hsm$) (\cite{Pr0} is a good exposition on skein modules). For the handful of cases where these two invariants fail, we compare the hyperbolic structures of the knot complements, which are harder to compute, but are complete invariants for hyperbolic knots (see Corollary \ref{cor:hyper}).

Let $M$ be an oriented 3-manifold and
let $\calLfr$ be the set of isotopy classes of unoriented framed links in $M$. 
Let $R = \Z[A^{\pm 1}]$ be the ring of Laurent polynomials in $A$ and 
let $R\calLfr$ be the free $R$-module generated by $\calLfr$.
Let $\calS$ be the ideal generated by the expressions:
\begin{align*}
\tag{Kauffman relator} \raisebox{-4pt}{\includegraphics[page=6]{images}} - A\;\raisebox{-4pt}{\includegraphics[page=7]{images}} -A^{-1}\;\raisebox{-4pt}{\includegraphics[page=8]{images}},\\
\tag{framing relator} L \sqcup \raisebox{-4pt}{\includegraphics[page=9]{images}} - (-A^2 - A^{-2}) L,
\end{align*}
where $\raisebox{-4pt}{\includegraphics[page=6]{images}}$, $\raisebox{-4pt}{\includegraphics[page=7]{images}}$, and $\raisebox{-4pt}{\includegraphics[page=8]{images}}$ are classes of links with representatives that are identical outside a small $3$-ball but look like the indicated diagrams inside it. Here blackboard framing is assumed. 

The Kauffman bracket skein module $\kbsm(M)$ is defined to be $R\calLfr$ modulo $\calS$. The $\kbsm$ is also called the second skein module and is often denoted as $\skbsm$~\cite{Pr0}.

We call a knot $K \subset M$ \df{affine} if it lies inside a 3-ball $B^3 \subset M$. 

\begin{prop}[Turaev~\cite{Tu}]\label{prop:KBSM-T}
$\kbsm(T)$ is freely generated by an infinite set of generators $\{x^n\}_{i=0}^\infty$, where $x^n$ is a parallel copy of $n$ longitudes of $T$ and $x^0$ is the affine unknot.
\end{prop}

\begin{prop}[Hoste, Przytycki~\cite{HP1}]\label{prop:KBSM-LPQ}
$\kbsm(\lpq)$ is freely generated by $\{x^n\}_{n=0}^{\lfloor p/2 \rfloor}$, where $x^n$ is a parallel copy of $n$ longitudes of $T \subset \lpq$ and $x^0$ is the affine unknot.
\end{prop}

If, for a given manifold $M$, the basis of $\kbsm(M)$ is known and $K$ is a knot in $M$, we denote by $\kbsm_M(K)$ the expression $K$ in $\kbsm(M)$ written in terms of the basis.

\begin{example}
	Take the knot $1_1$ in the solid torus from Appendix A and resolve the crossing using the Kauffman relator in $\kbsm(T)$:
	$$ \bigg[ \raisebox{-14pt}{\includegraphics[page=28]{images}} \bigg] = \
	A \; \bigg[ \raisebox{-14pt}{\includegraphics[page=30]{images}} \bigg] + \
	A^{-1}\bigg[ \raisebox{-14pt}{\includegraphics[page=29]{images}} \bigg]. $$
	
	Since we end up with an expression in the basis we can write $\kbsm_T(1_1) = A\,x^2 + A^{-1}\,x^0.$
\end{example}



In contrast to the $\kbsm$, see for example~\cite{Tu, HP1, Pr1, mmm, Mr, MM}, the $\hsm$ is not a widely studied knot invariant. Uncoincidentally, the $\hsm$ is a much stronger invariant and is much more difficult to compute, see for example \cite{Tu, HK, La1, BM4}.
The $\hsm$ has been calculated for the solid torus using diagrammatic methods \cite{Mr} as well as algebraic methods \cite{Tu,HK, La1, DL2}, for $S^1 \times S^2$~\cite{HP2, Mr}, and recently for the family of lens spaces of type $L(p,1)$~\cite{BM4}, see also~\cite{DLP}. The computation of $\hsm(L(p,q)), q > 1$ is still an open problem, but evidence suggests that it is free with the same basis as $\hsm(L(p,1))$.

Let $M$ again be an oriented 3-manifold and
let $\calLor$ be the set of isotopy classes of oriented links in $M$ to which we also add the empty knot $\emptyset$.
Let $R = \Z[v^{\pm 1}, z^{\pm 1}]$ be the ring of Laurent polynomials in variables $v$ and $z$, and 
let $R\calLor$ be the free $R$-module generated by $\calLor$.
Let $\calS$  be the ideal generated by the expression:
\begin{align*}
\tag{HOMFLYPT relator} v^{-1}\;\raisebox{-4pt}{\includegraphics[page=10]{images}} - v\;\raisebox{-4pt}{\includegraphics[page=11]{images}} - z\;\raisebox{-4pt}{\includegraphics[page=12]{images}}.
\end{align*}
We also add to $\calS$ the expression involving the empty knot: 
\begin{align*}
\tag{empty knot relator} v^{-1}\emptyset - v\emptyset - z\,\raisebox{-4pt}{\includegraphics[page=13]{images}}.
\end{align*}
The \df{HOMFLYPT skein module} $\hsm(M)$ is $R\calLor$ modulo $\calS$.
The $\hsm$ is also called the third skein module and is often denoted by $\shsm$~\cite{Pr0}.

\begin{prop}[Turaev~\cite{Tu}]\label{prop:HSM-T}
$\hsm(T)$ is freely generated by an infinite set of generators $\mathcal{B}=\{t^{i_1}_{k_1} \ldots t^{i_s}_{k_s} \, | \, s\in\N,
k_1< \cdots <k_s\in\Z \setminus \{0 \},i_1, \ldots, i_s\in \N\}\cup\{\emptyset\}$, where
$\emptyset$ is the empty knot. For $k>0$, $t_k$ is the oriented knot in $T$ 
representing $k$ in $\pi _1(T)\cong\Z$ that has an ascending diagram with $k-1$ positive crossings. For $k<0$, $t_k$ is $t_{|k|}$ with reversed orientation.
\end{prop}

\begin{figure}[htb]
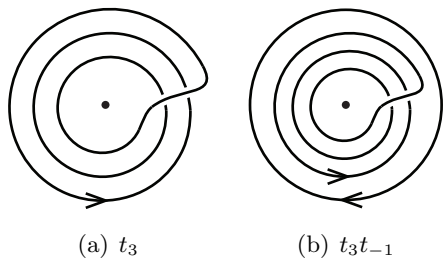

\centering
\subfigure[$t_3$]{\begin{overpic}[page=24]{images}\end{overpic}}
\subfigure[$t_3 t_{-1}$]{\begin{overpic}[page=25]{images}\end{overpic}} 
\caption{Two generators of $\hsm(T)$.}
\label{fig:HOMFLYPT-torus}
\end{figure}

For example, $t_3$ and $t_{-1}t_3$ are presented in Figure~\ref{fig:HOMFLYPT-torus}.

\begin{prop}[Gabrov\v sek, Mroczkowski~\cite{BM4}]\label{prop:HSM-LP1}
$\hsm(\lpo)$ is freely generated by an infinite set of generators $\calB_{p}=\{t^{i_1}_{k_1} \ldots t^{i_s}_{k_s}:s\in\N,
k_1< \cdots <k_s\in\Z \setminus \{ 0\},-\frac{p}{2}<k_1< \cdots <k_s\le \frac{p}{2},i_1, \ldots, i_s\in 
\N\}\cup\{\emptyset\}$, where
$\emptyset$ is the empty knot and $t_k$ are knots with diagrams equal to those in Proposition~\ref{prop:HSM-T}.
\end{prop}




Having fixed a basis of $\hsm(M)$ and having a knot $K \subset M$  we denote by $\hsm_M(K)$ the expression of $K$ in $\hsm(M)$ written in terms of the basis. 

\begin{prop}\label{prop:chir1}
Let $K$ be a knot in $T$ and $\overline{K}$ its mirror. Then the evaluations $\kbsm_T(K)$ and $\kbsm_T(\overline{K})$ are related by a change of variable $A \leftrightarrow A^{-1}$.
\end{prop}

\begin{prop}\label{prop:chir2}
	Let $K$ be a knot in $T$ and $\overline{K}$ its mirror. The evaluations $\hsm_T(K)$ and $\hsm_T(\overline{K})$ are related by a change of variable $v \leftrightarrow -v^{-1}$.	
\end{prop}

The first proposition follows directly from the Kauffman relator and the second proposition follows directly from the HOMFLYPT relators.

\subsection{Comparing the hyperbolic structures}

As we will see in Section~\ref{sec:classT} not all knots are distinguished by their skein module evaluations. For these knots we compare the complements with an approach similar to \cite{HTW}, but we face greater computational difficulties since the ambient spaces are lens spaces. 

A knot $K$ in a 3-manifold $M$ is called \df{hyperbolic} if the complement $M-K$ admits a complete Riemannian metric of constant Gaussian curvature $-1$. Two hyperbolic knots are ambient isotopic if and only if the hyperbolic structures of their complements are the same. This follows from the following theorems.

\begin{thm}[Gordon-Luecke~\cite{gl}]  \label{thm:gl}¸
	Two knots are isotopic if and only if their complements are orientation-preserving homeomorphic.
\end{thm}

\begin{thm}[Mostow rigidity theorem~\cite{pra}] \label{thm:pra}
	Let $M$ and $N$ be complete finite-volume hyperbolic $n$-manifolds of dimension $n > 2$. If there exists an isomorphism $f:\pi_1(M) \rightarrow \pi_2(M)$ then it is induced by a unique isometry from $M$ to $N$.
\end{thm}

\begin{thm}[Epstein-Penner-Weeks~\cite{ep, W}] \label{thm:epw}
There exists a canonical cell decomposition of a cusped hyperbolic $n$-manifold.
\end{thm}
	
\begin{cor}[\cite{HTW}]\label{cor:hyper}
Two hyperbolic knots $K_1$ and $K_2$ in a 3-manifold $M$ are isotopic if and only if there is an isomorphism from the canonical triangulation of $M-K_1$ to $M-K_2$.
\begin{proof}
It follows from Theorem~\ref{thm:pra} that if a knot is hyperbolic, then the Riemannian metric is unique. Combined with Theorem~\ref{thm:gl} it follows that two knots are isotopic if and only if their complements are isometric. By Theorem~\ref{thm:epw} it follows that a knot's complement has a unique triangulation, which arises from the unique metric.
\end{proof}
\end{cor}

Fortunately, the computational tool SnapPy~\cite{snappy} is able to verify if a knot is hyperbolic, triangulate the complement and combinatorically find isomorhisms between two triangulations.
 
\section{Gauss codes}\label{sec:gauss}

We use Gauss codes as the computer data structure to store knots. A \df{Gauss word} of length $n$ is a word on the alphabet $\{ \pm 1, \pm 2 , \ldots \pm n \}$, where each letter appears exactly once.

Let $D$ be a diagram of a knot in $S^3$ with $n$ crossings. The Gauss word of length $n$ associated with $D$ is obtained by the following steps~\cite{G}:
\begin{enumerate}
\item Enumerate the crossings of $D$ from $1$ to $n$.
\item Orient $D$ (if not already oriented) and choose an initial point on it.
\item Starting with an empty word, travel from the initial point back to it according to the orientation, appending the letter $k$ to the end of the word when passing through an overcrossing and the letter $-k$ when passing through an undercrossing.	
\end{enumerate}

\begin{figure}[htb]
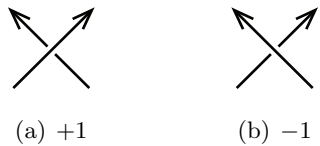

\centering
\subfigure[$+1$]{\hspace{0,7cm}\begin{overpic}[page=33]{images}\end{overpic}\hspace{0,7cm}}\label{fig:sp}
\subfigure[$-1$]{\hspace{0,7cm}\begin{overpic}[page=34]{images}\end{overpic}\hspace{0,7cm}}\label{fig:sn}
\caption{The crossing sign.}
\label{fig:sgn}
\end{figure}

In general we can reconstruct a knot from its associated Gauss word only up to its mirror image, but having the knowledge of the knot's \df{crossing signs} enables us to reconstruct the knot completely.
Each crossing of a diagram $D$ can be assigned its corresponding crossing sign according to the right-hand rule illustrated in Figure~\ref{fig:sgn}. 

\begin{figure}[htb]
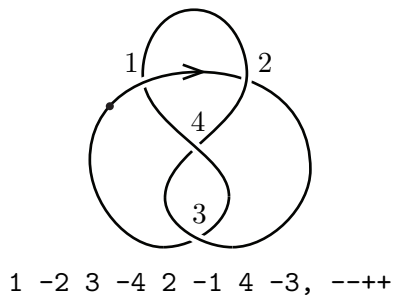

\centering
\begin{overpic}[page=26]{images}\put(17,72){1}\put(69,72){2}\put(43.5,13){3}\put(43,49){4}\end{overpic}
\gauss{1 -2 3 -4 2 -1 4 -3, --++}
\caption{The oriented figure eight knot and its Gauss code.}
\label{fig:f8}
\end{figure}

The \df{Gauss code} of a diagram $D$ is the Gauss word of $D$ followed by the sequence of crossings signs in the order of the enumeration (Figure~\ref{fig:f8}).

A Gauss code that does not realize a knot is called \df{nonrealizable}~\cite{Marx}.

We will now recall from~\cite{BM} how to adopt the notation of a Gauss code for $T$.
To extend the Gauss code to a code of a knot in $T$, we take a punctured disk diagram $D$ and keep track of the $0$- and $\infty$-regions.
An \df{extended Gauss code} is a Gauss code followed by two sequences: the arcs bounding the $0$-region and the arcs bounding the $\infty$-region, where we enumerate the arcs in the direction of orientation, starting with arc $0$ positioned right after the first crossing when crossing from the over-arc (see Figure~\ref{fig:pddex}).
We will omit the word ``extended'' if the variant of the Gauss code will be evident from the context.

\begin{figure}[htb]
\centering
\begin{overpic}[page=27]{images}\put(17,72){1}\put(69,72){2}\put(43.5,13){3}\put(42,29){4}
\put(50,75){\scriptsize 0}\put(60,50){\scriptsize 3}\put(92,32){\scriptsize 1}
\put(27.5,47){\scriptsize 5}\put(-2,34){\scriptsize 7}\put(42,99){\scriptsize 4}
\put(27,20){\scriptsize 2}\put(60,20){\scriptsize 6}
\end{overpic}
\gauss{1 -2 3 -4 2 -1 4 -3, --++, 0 3 5, 1 4 7}
\caption{A punctured figure eight knot and its corresponding extended Gauss code.}
\label{fig:pddex}
\end{figure}

The \df{length} of a Gauss code is the length of the underlying Gauss word (i.e. the number of crossings of the corresponding diagram).  

We introduce a total order on the set of (extended) Gauss codes. The order is the lexicographical order obtained by considering in turn:
\begin{enumerate}
\item the length of the Gauss code $n$,
\item the lexicographical ordering $1<-1<2<-2<\ldots<n<-n$  of the Gauss word,
\item the lexicographical ordering $+ < -$ of signs,
\item the lexicographical ordering $0<1<2<\ldots<2n-1$ of the arcs bounding the $0$-region,
\item the lexicographical ordering $0<1<2<\ldots<2n-1$ of the arcs bounding the $\infty$-region.
\end{enumerate}

If $W$ and $W'$ are two Gauss codes that represent the same knot and $W' < W$, we say that $W'$ is a \df{reduction} of $W$.
A Gauss code that allows a reduction is called \df{reducible}.

For example, the Gauss code calculated in Example~\ref{fig:pddex} is reducible, since,
if we choose the initial point on arc $1$, the Gauss code of the diagram becomes
\gauss{1 -2 3 -4 2 -1 4 -3, ++--, 1 3 6, 2 5 7.}

\section{The classification algorithm}\label{sec:alg}

We describe the algorithm used to classify knots in the solid torus and lens spaces, which will provide us with the desired knot tables.

We only classify non-affine knots, since affine knots agree with those in $S^3$ (by the fact that the knot theory of $S^3$ embeds in the Knot theory of $T$) and have already been classified by classical knot tables.

We recall from~\cite{BM} that unlike the classical case, there is no well defined connected sum operation for knots in $T$. But for two oriented knots, where at least one of them is affine, the connected sum operation is well defined. The same holds for $L(p,q)$.
A \df{prime knot} in $T$ or $L(p,q)$ is a knot that cannot be expressed as a connected sum of two non-trivial knots, where at least one of them is affine.
A well defined connected sum operation implies that the connected sum decomposition is unique (cf.~\cite{Schubert}). It is therefore enough to classify only prime knots since all other can be obtained by composing them with affine knots. 

\subsection{Classifying knots in the solid torus}\label{sec:classT}

In order to classify knots in the solid torus $T$ with at most $n$ crossings, we use the following steps:

{\bf Step 1. Find all realizable Gauss codes.}
A Gauss code is in its \df{canonical form} if it cannot be reduced by any combination of the following operations:
\begin{itemize}
\item renaming the letters of the Gauss code (this corresponds to renaming the crossings in $D$),
\item a cyclic shift of the Gauss code (this corresponds to choosing an initial point in $D$),
\item reversing the Gauss code (this corresponds to reversing the orientation of $D$).
\end{itemize}
Note that a knot diagram can be uniquely represented by the canonical form of its Gauss code.

We start off by generating all Gauss codes in the lexicographical order.
At this step we also eliminate Gauss codes (diagrams) that:
\begin{itemize}
	\item are not in their canonical form (their canonical partner already exists lower in the list),
	\item are nonrealizable (which is easy to check with methods described in~\cite{L, Marx}), 
	\item the $0$-region and the $\infty$-region coincide (affine knot),
	\item the $0$-region and the $\infty$-region are adjacent (the sum of an affine knot and the non-affine unknot). An exception is the non-affine unknot which we must count.
\end{itemize}


The number of Gauss codes up to a length of 7 is presented in Table~\ref{table:ngc}: the first row of the table is given by the formula $N=(2n)!\cdot 2^n$, other rows are determined experimentally.

\begin{table}[h]
\caption{Number of Gauss codes up to length 7.} \label{table:ngc}
\centering
\begin{tabular}{|l|c|c|c|c|c|c|c|}
\hline
{\bf \it n} & 1 & 2 & 3 & 4 & 5 & 6  & 7 \\
\Xhline{2\arrayrulewidth}
Gauss codes                            & $4$ & $96$   & $5760$ &$6.5 \!\cdot\! 10^5$ & $1.2 \!\cdot\! 10^8$&$3.1 \!\cdot\! 10^{10}$ & $1.1 \!\cdot\! 10^{13}$\\
Realizable                                & $4$& $24$   & $432$   & $13344$   & $7.1 \!\cdot\!10^5$            & $5.2 \!\cdot\! 10^7$ & $4.9 \!\cdot\! 10^{9}$\\
Canonical real.                  &$2$ & $6$    & $36$    & $278$       & $3.0 \!\cdot\! 10^3$                &$3.6 \!\cdot\! 10^4$ & $4.8 \!\cdot\! 10^{5}$\\
Extended real.                  & $36$& $384$& $10800$&$4.8 \!\cdot\! 10^5$ & $3.5 \!\cdot\! 10^7$&$3.3 \!\cdot\! 10^9$ & $4.0 \!\cdot\! 10^{11}$\\
Can. ext. real. & $10$& $68$&$714$     &$9392$      &$1.4 \!\cdot\! 10^5$ &            $2.3 \!\cdot\! 10^6$ & $4.0 \!\cdot\! 10^{7}$\\
\hline
\end{tabular}
\end{table}

\textbf{Step 2. Partitioning the knots that share the same $\kbsm$.} Let $\gc$ be the set of all Gauss codes calculated in Step 1. We partition $\gc$ into partitions $P_0, P_1,\ldots$ 
so that all codes in each partition are evaluated the same in $\kbsm(T)$:
$$\kbsm(W_1) = \kbsm(W_2) \Leftrightarrow W_1 \in P_i \wedge W_2 \in P_i \mbox{ for } W_{1,2}\in \gc.$$

The first knot (with respect to the ordering) in a partition cannot be reduced, since we would get a lower order knot in the partition, which is a contradiction.
We therefore mark the first knot in each partition as prime.




\textbf{Step 3. Finding isotopic knot diagrams.} 

In this step we perform a breadth-first search (BFS) algorithm to systematically perform Reidemeister moves on all Gauss codes not marked as prime in the previous step.

If in this process a Gauss word $W$ reduces to a Gauss code $W' < W$, we can eliminate it from the set of possible prime knots, since the reduced code $W'$ already exists lower in the lexicographical order. The codes are stored in binary search trees (BSTs). We use BSTs since they realise sets (i.e. two diagrams will never be processed simultaneously) and time complexities of operations lie in our favour: insertation of an element takes $O(\ln n)$ time, retrieving elements in order and finding the minimal element take constant time $O(1)$.

In this step we do not include only Reidemeister moves as the set of allowed moves, but also other isotopies: flypes, meridional rotations and ``moving a strand through several crossings''. The first two isotopies are relatively easy to perform on a Gauss code, but would otherwise involve a large amount of Reidemeister moves. The third isotopy is a result of performing several consecutive $\riii$ moves whenever we can.

After the process we eliminate connected sums. Since it holds that $$\kbsm(K_1 \# K_2) = \kbsm(K_1) \cdot \kbsm(K_2),$$ we check for each prime knot candidate if its $\kbsm$ evaluation factors. It turns out that all such ``$\kbsm$-divisible'' knots are connected sums, which we check by hand.

Since there is a well defined mirror operation in $T$ (see Section~\ref{sec:kdasm}), we also establish which knots in $T$ are amphichiral.
Detecting amphichirality is easy: a knot $K$ is amphichiral if its mirror $\overline{K}$ belongs to the same partition as $K$.

The above process works for all knots up to four crossings. For five crossings there are a few partitions left with two elements. We deal with these knots in Section~\ref{sec:results}.

\subsection{Classifying knots in \boldmath{$\lpq$}}\label{sec:classLPQ}

As seen in Section~\ref{sec:lpq}, a punctured disk diagram of a knot in $\lpq$ can be thought of as a punctured disk diagram in $T$ accompanied by the additional $\slajd$ move.

The algorithm for classifying knots in $\lpq$ is therefore the same as the one for $T$, except that we add the $\slajd$ move to the set of Reidemeister moves and we only search for prime knots within the set of prime knots of $T$.

Again the algorithm works only for knots up to four crossings. We deal with the 5-crossing exceptions in the following section.

\begin{remark} There are prime knots in $T$ with more than 5 crossings which reduce to 5 or less crossings in $L(p,q)$, but the latter knot has a diagram in $L(p,q)$ with up to 5 crossings and has therefore already been acocunted for. Similarly, if there were reducible, affine, or composite knots in $T$, which are prime in $L(p,q)$ up to 5 crossings, we would not have to include them in our search, since their reduced diagram is already in the table.
\end{remark}






\section{The results}\label{sec:results}

\subsection{The solid torus}\label{sec:resT}

The number of non-affine prime knots in the solid torus up to 5 crossings is presented in Table~\ref{table:ctorus}.
The knot table is presented in Appendix~A.

\begin{table}[H]
\caption{The number of non-affine prime knots in the solid torus.} \label{table:ctorus}
\centering
\begin{tabular}{| c | l l |} 
\hline
$n$ &  \multicolumn{2}{l |}{Number of prime knots } \\
\hline
\;\;{\bf 0}\;\; & 1 & (0 mirror pairs + 1 amphichiral)\\
\;\;{\bf 1}\;\; & 2 & (1 mirror pair + 0 amphichiral)\\
\;\;{\bf 2}\;\; & 5 & (2 mirror pairs + 1 amphichiral)\\
\;\;{\bf 3}\;\; & 16 & (8 mirror pairs + 0 amphichiral)\\
\;\;{\bf 4}\;\; & 50 & (23 mirror pairs + 4 amphichiral)\\
\;\;{\bf 5}\;\; & 190 & (95 mirror pairs + 0 amphichiral) \\
\hline
\end{tabular}
\end{table}
Up to five crossings there are six amphichiral knots: $0_1$, $2_2$, $4_{13}$, $4_{14}$, $4_{21}$, and $4_{27}$. We pose the question if there exists a reduced amphichiral knot with an odd number of crossings.

There are five pairs of knots that are not distinguished by the $\kbsm$, $4_{23}$, $\overline{5_1}$; $5_1$, $\overline{4_{23}}$; $5_{26}$, $5_{27}$; $\overline{5_{26}}$, $\overline{5_{27}}$; and $5_{76}$, $\overline{5_{76}}$:
\begin{alignat*}{4}
&\kbsm_T(4_{23})  &=& \kbsm_T(\overline{5_1})  &= &\,  ( - A^{12} + A^{8} )\,x^{3} + ( 2A^{12} - A^{8} )\, x, \\
&\kbsm_T(5_1)  &=& \kbsm_T(\overline{4_{23}})  &= &\,  ( A^{-8} - A^{-12} ) x^{3} + ( - A^{-8} + 2A^{-12} )\, x, \\
&\kbsm_T(5_{26})  &=& \kbsm_T(5_{27})  &=&\,  ( A^{-4}-2A^{-8} + 2A^{-12} - A^{-16} )\, x^{3}, \\
& & & & & + ( - 1 + A^{-4} + A^{-8} - A^{-16} + A^{-20} )\, x, \\
&\kbsm_T(\overline{5_{26}})  &=& \kbsm_T(\overline{5_{27}})  &=&\, ( A^{-4} - 2A^{-8} + 2A^{-12} - A^{-16} )\, x^{3}, \\ 
& & & & & + ( - 1 + A^{-4} + A^{-8} - A^{-16} + A^{-20} )\, x \\
& \kbsm_T(5_{76})  &=& \kbsm_T(\overline{5_{76}})  &= &\,  ( A^{4} - 1 + A^{-4} )\, x^{3} + ( - 2A^{4} + 1 - 2A^{-4} )\, x.
\end{alignat*}
Three of these pairs are distinguished by the $\hsm$:
\begin{alignat*}{2}
 &\hsm_T(4_{23})  &=&  {-\tfrac{z}{v^3}}\,t_{-1}t_{2} + (\tfrac{z^2}{v^2} + \tfrac{1}{v^2})\,t_{1}, \\
 &\hsm_T(\overline{5_1})  &=& (\tfrac{z^4}{v^4} + \tfrac{2z^2}{v^4})\,t_{-1}t_{1}^{2} - (\tfrac{z^3}{v^5} + \tfrac{z}{v^5})\,t_{-1}t_{2} + (\tfrac{z^2}{v^4} + \tfrac{2}{v^4} - \tfrac{1}{v^6})\,t_{1}, \\
 &\hsm_T(5_1)  &=& v^{4}z^{2}\,t_{-1}t_{1}^{2} + (v^{3}z^{3} + v^{3}z)\,t_{-1}t_{2} +(- v^{6} + v^{4}z^{2} + 2v^{4})\,t_{1}, \\
 &\hsm_T(\overline{4_{23}})  &=& {-v^{2}z^{2}}\,t_{-1}t_{1}^{2} + vz\,t_{-1}t_{2} + (v^{2}z^{2} + v^{2})\,t_{1}, \\
 &\hsm_T(5_{76}) &=& {-vz\,t_{1}t_{2}} + (- z^{2} + \tfrac{1}{v^2})\,t_{3} \\
 &\hsm_T(\overline{5_{76}}) &=&  (z^{2} - \tfrac{z^4}{v^2} - \tfrac{z^2}{v^2})\,t_{1}^{3} +(- \tfrac{2z}{v} + \tfrac{2z^3}{v^3} + \tfrac{z}{v^3})\,t_{1}t_{2} + (\tfrac{1}{v^2} - \tfrac{z^2}{v^4})\,t_{3}.
\end{alignat*}
It follows from Proposition~\ref{prop:chir1} that if a knot is amphichiral, the $\kbsm$ evaluation is symmetric in variable $A$. Similarly, it follows from Proposition~\ref{prop:chir2} that if a knot is amphichiral, the $\hsm$ evaluation is symmetric in variable $v$ up to sign. Note that for the knot $5_{76}$ the $\kbsm$ suggests amphichirality, but the $\hsm$ shows otherwise.

The pair $5_{26}, 5_{27}$ and its mirrors share the same $\hsm$ evaluation:
\begin{alignat*}{3}
&\hsm_T(5_{26})  &=& \hsm_T(5_{27})  &=&  (-v^{2}z^{4} - v^{2}z^{2})\,t_{-1}t_{1}^{2} \\
& & & & & + (vz^{3} + vz)\,t_{-1}t_{2} + (v^{4}z^{2} + v^{2})\,t_{1}, \\
&\hsm_T(\overline{5_{26}})  &=& \hsm_T(\overline{5_{27}}) &=& (- \tfrac{z^3}{v^3} - \tfrac{z}{v^3})\,t_{-1}t_{2} + (\tfrac{1}{v^2} + \tfrac{z^2}{v^4})\,t_{1}. 
\end{alignat*}
These pairs are hyperbolic and there is no isomorphism between their canonical triangulations (see \cite{code}).

\begin{figure}[ht]
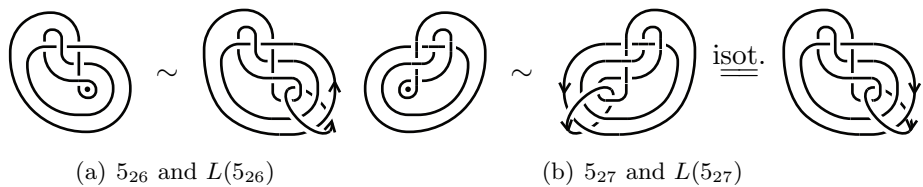

\centering
\subfigure[$5_{26}$ and $L(5_{26})$]{\begin{overpic}[page=31]{images}\end{overpic}  \raisebox{0.85cm}{$\;\sim\;$}  \begin{overpic}[page=32]{images}\end{overpic}}
\subfigure[$5_{27}$ and $L(5_{27})$]{\begin{overpic}[page=41]{images}\end{overpic}  \raisebox{0.85cm}{$\;\sim\;$}  \begin{overpic}[page=42]{images}\end{overpic}   \raisebox{0.85cm}{$\stackrel{\mbox{isot.}}{=\!=}$}  \begin{overpic}[page=43]{images}\end{overpic}}
\caption{$5_{26}$ and $5_{27}$ and their links.}
\label{fig:526}
\end{figure}

The pair $5_{26}$ and $5_{27}$ (and their mirrors) are of particular interest. Observe from Figure~\ref{fig:526} that $5_{26}$ and $5_{27}$ are flips of each other (note that in Figure~\ref{fig:526} the knot $5_{27}$ differs from the one in Appendix A by a meridional rotation).
An oriented knot $K$ in the solid torus naturally corresponds to a two-component oriented link $L(K) \subset S^3$ with one trivial component, the extra trivial component being the meridian of $\rob T$~\cite{HK, LR1, La3}. The knots $5_{26}$ and $5_{27}$ and their corresponding links are presented in Figure~\ref{fig:526}. Note that we have oriented $5_{26}$ and $5_{27}$ in such a way that the linking numbers match, $lk(L(5_{26}))=lk(L(5_{27}))$. If we rotate $L(5_{27})$ by $\pi$ in $\R^3$ we see that $L(5_{26}) = -L(5_{27})$. The question if $5_{26} = 5_{27}$ is thus equivalent to the question if $L(5_{26}) = -L(5_{26})$. 
A link with the property that reversing the orientations on both components yields the same link is called an \df{invertible link}. 
Non-invertibility is a very difficult question in knot theory and only few examples of such links have been found so far (see Whitten~\cite{whitten, whitten2}).
Whitten used algebraic methods to prove the existance of non-invertible links. He found examples of links where there is no automorphism of the knot group $\pi_1(S^3 - K_1)$
that takes $[K_2] \in \pi_1(S^3 - K_1)$ to $[-K_2] \in \pi_1(S^3 - K_1)$ and takes the meridian of $K_1$, which is a generator of $\pi_1(S^3 - K_1)$, to its inverse. Whitten concluded that if the link cannot be inverted algebraically, it cannot be inverted by isotopy.

For the link $L(5_{26}) = K \cup U$, where $U$ is the trivial component and $K$ is the knotted component, it holds that
$\pi_1(S^3 - U) = \langle m \rangle \cong \Z$ and $[K] = m \in \langle m \rangle$. There exists an obvious automorphism $\langle m \rangle \rightarrow \langle m \rangle$, namely the one that takes $m \mapsto -m$ and also takes $[K] \mapsto [-K]$. The link can therefore be algebraically inverted and we cannot conclude by this method that it is non-invertible.

\subsection{Lens spaces}\label{sec:resLPQ}

We tabulate non-affine prime knots in the lens spaces $\lpq$, $p \leq 12$ and argue that for $p > 12$ the table agrees with that of the solid torus at least for the knots distinguished either by the $\kbsm$ or $\hsm$. Up to homeomorphism there are exactly 19 lens spaces with $p \leq 12$:
$L(2,1)$, $L(3,1)$, $L(4,1)$, $L(5,1)$, $L(5,2)$, $L(6,1)$, $L(7,1)$, $L(7,2)$, $L(8,1)$, $L(8,3)$, $L(9,1)$, $L(9,2)$, $L(10,1)$, $L(10,3)$, $L(11,1)$, $L(11,2)$, $L(11,3)$,  $L(12,1)$, and $L(12,5)$~\cite{R1}.

The \df{winding number} $\wind(K)$ of a knot $K \subset T$ is the integer $[K] \in \pi_1(T) \cong \Z$, or equivalently, the integer $[D] \in \pi_1(\R^2 \setminus \{ \cdot \} ) \cong \Z$ where $D$ is a diagram of $K$ lying in the punctured plane. Here $[K]$ (resp. $[D]$) represents the homotopy class of $K$ (resp. $D$).

\begin{prop}\label{prop:kbsm}
Let $i: T \hookrightarrow \lpq$ be the standard inclusion of the solid torus $T$ in the lens space $\lpq$.
If two knots $K_1$ and $K_2$ with at most $n$ crossings have distinct $\kbsm$s in $T$, then the knots $i(K_1)$ and $i(K_2)$ have distinct $\kbsm$s in $\lpq$, $p \geq 2(n+1)$.
\begin{proof}
Since a knot $K \subset T$ that allows a diagram $D$ with $n$ crossings has the winding number bounded by $| \wind(K) | \leq n+1$ ($D$ can go up to $(n+1)$ times around the puncture and thus needs $n$ crossings to complete the cycle), it holds that $\kbsm_T(K) = \sum_{i=0}^{n+1} A_i\,x^i$, for some $A_i \in R$.
Since $p \geq 2(n+1) \Rightarrow \lfloor p/2 \rfloor \geq n+1$ and the generators $\{\hat{x}^n\}_{n=0}^{\lfloor p/2 \rfloor}$ of $\kbsm(\lpq)$ are induced by the inclusion, $i_*(x^i) = \hat{x}^i, 0 \leq i \leq \lfloor p/2 \rfloor$, it holds that $\kbsm_{\lpq}(i(K)) = \sum_{i=0}^{n+1} A_i\,\hat{x}^i$, $A_i \in R$. That is, the knots $K$ and $i(K)$, although lying in different spaces, have an equal $\kbsm$ expression.
\end{proof}
\end{prop}

A similar proposition can be made for the $\hsm$ of the spaces $L(p,1)$.

\begin{prop}
Let $i: T \hookrightarrow \lpo$ be the standard inclusion of the solid torus $T$ in the lens space $\lpo$.
If two knots $K_1$ and $K_2$ with at most $n$ crossings have distinct $\hsm$s in $T$, then the knots $i(K_1)$ and $i(K_2)$ have distinct $\hsm$s in $\lpo$, $p > 2(n+1)$.
\begin{proof}
Since the generators $\{t^{i_1}_{k_1}\ldots t^{i_s}_{k_s}:s\in\N,
k_1<..<k_s\in\Z \setminus \{ 0\},-\frac{p}{2}<k_1<..<k_s\le \frac{p}{2},i_1..i_s\in 
\N\}\cup\{\emptyset\}$ of $\hsm(\lpo)$ are induced by the inclusion of the generators of $\hsm(T)$, by the same arguments used in Proposition~\ref{prop:kbsm} and by the fact that $p>2(n+1) \Rightarrow \frac{p}{2} > n$, we conclude that the expression $\hsm_T(K)$ is equal to the expression $\hsm_{\lpo}(i(K))$.
\end{proof}
\end{prop}


The number of non-affine prime knots up to five crossings in $\lpq$ is presented in Table~\ref{table:numknotstable}. The subset of knots from the knot table of the solid torus in Appendix~A that are either reducible, affine, or composites in $\lpq$ are tabulated in Appendix~B. Note that we also tabulate the mirrors, since there is no well-defined mirror operation in $\lpq$. Let us summarize the table in Appendix~B:
\begin{itemize}
	\item in $L(2,1)$ there are 47 knots that reduce to prime knots, in addition, the knots $1_1$ and $5_{74}$ reduce to the affine unknot (in Appendix A we can see directly that $1_1$ is the result of the slide move on the unknot), $5_{24}$ is the affine trefoil, and $5_{77}$ is the connected sum of the non-affine unknot and the trefoil,
	
	\item in $L(3, 1)$ there are 38 knots that reduce to prime knots and the knots $2_3$ and $5_{71}$ reduce to the affine unknot (again it can be directly seen that $2_3$ is the result of the slide move on the unknot),
	
	\item in $L(4, 1)$ there are 21 knots that reduce to prime knots and the knot $3_6$ is the affine unknot,
	
	\item in $L(5, 1)$: there are 12 knots that reduce to prime knots and the knot $4_{26}$ is the affine unknot,
	
	\item in $L(5, 2)$: there are 4 knots that reduce to prime knots,
	
	\item in $L(6, 1)$: there are 8 knots that reduce to prime knots and the knot $5_{85}$ si the affine unknot,
	
	\item in higher lens spaces there are only a handful of reducible knots, if any.

\end{itemize}

\begin{table}[ht]
\caption{Number of non-affine prime knots in $\lpq$ with up ot 5 crossings.} \label{table:numknotstable}
\centering
\begin{tabular}{| c | c c c c c c c |} 
\hline
$n$ & $L(2,1)$ & $L(3,1)$ & $L(4,1)$ & $L(5,1)$ & $L(5,2)$ & $L(6,1)$  & $L(7,1)$ \\ \hline
{\bf 0} & 1 & 1 & 1 & 1 & 1 & 1 & 1 \\
{\bf 1} & 1 & 1 & 2 & 2 & 2 & 2 & 2 \\
{\bf 2} & 4 & 4 & 4 & 4 & 4 & 5 & 5 \\
{\bf 3} & 12 & 13 & 14 & 15 & 16 & 15 & 15 \\
{\bf 4} & 40 & 42 & 45 & 46 & 49 & 48 & 49 \\
{\bf 5} & 155 & 163 & 176 & 183 & 188 & 184 & 186 \\
\specialrule{1.5pt}{1pt}{1pt}
$n$ & $L(7,2)$ & $L(8,1)$ & $L(8,3)$ & $L(9,1)$ & $L(9,2)$ & $L(10,1)$ & $L(10,3)$ \\ \hline
{\bf 0} & 1 & 1 & 1 & 1 & 1 & 1 & 1 \\
{\bf 1} & 2 & 2 & 2 & 2 & 2 & 2 & 2 \\
{\bf 2} & 5 & 5 & 5 & 5 & 5 & 5 & 5 \\
{\bf 3} & 15 & 16 & 16 & 16 & 16 & 16 & 16 \\
{\bf 4} & 50 & 49 & 50 & 49 & 49 & 50 & 50 \\
{\bf 5} & 189 & 188 & 190 & 189 & 190 & 189 & 190 \\
\specialrule{1.5pt}{1pt}{1pt}
$n$ & $L(11,1)$ & $L(11,2)$ & $L(11,3)$ & $L(12,1)$ & $L(12,5)$ & \multicolumn{2}{c |}{$\lpq,\, p \geq 13$}  \\ \hline
{\bf 0} & 1 & 1 & 1 & 1 & 1 & \multicolumn{2}{c |}{1} \\
{\bf 1} & 2 & 2 & 2 & 2 & 2 & \multicolumn{2}{c |}{2} \\
{\bf 2} & 5 & 5 & 5 & 5 & 5 & \multicolumn{2}{c |}{5} \\
{\bf 3} & 16 & 16 & 16 & 16 & 16 & \multicolumn{2}{c |}{16} \\
{\bf 4} & 50 & 50 & 50 & 50 & 50 & \multicolumn{2}{c |}{50} \\
{\bf 5} & 189 & 189 & 190 & 190 & 190 & \multicolumn{2}{c |}{$190^{\ast}$} \\
\hline
\end{tabular}
\\[0.2cm]
$^\ast$ Conjecture
\end{table}

Since the $\kbsm$s are unique for all knots up to 4 crossings in $T$ and all knots up to four crossings in $\lpq$, $p < 10$, we conclude by Proposition~\ref{prop:kbsm} that the $\kbsm$ distinguishes these knots in all lens spaces $\lpq$.


In Table~\ref{table:sharebracket} we present groups of inequivalent knots that are not distinguished by the $\kbsm$.
For $L(p,1)$ all but the pairs $5_{26}$, $5_{27}$ and $\overline{5_{26}}$, $\overline{5_{27}}$ are distinguished by the $\hsm$.

For $(p,1)$, where $0 < p \leq 12$, the pairs $5_{26}$, $5_{27}$ and $\overline{5_{26}}$, $\overline{5_{27}}$ are distinguished by their hyperbolic structure (recall Corollary~\ref{cor:hyper}). We verify this by taking the mixed link complements $S^3 - L(5_{26})$ and $S^3 - L(5_{27})$, fill the cusp of the trivial component by a $1/p$-Dehn filling, verify the complements are hyperbolic, calculate the canonical triangulation of the manifold and combinatorically compute the isomorphisms between the two triangulations. The computations show that for $0 < p \leq 12$ there are no isomorphisms between the triangulations from where we conclude that the pairs $5_{26}$, $5_{27}$ and $\overline{5_{26}}$, $\overline{5_{27}}$ are inequivalent. All of the above operations are performed with built-in SnapPy methods, the computer code with full outputs is available at~\cite{code}.
For $p\geq 13$ we cannot conclude that they represent different knots since our computational methods do not work for general $p,q$-s.

For $q>1$ and $0 < p \leq 12$ all pairs from Table~\ref{table:sharebracket} are distinguished by their hyperbolic structure, with same methods as in the previous paragraph, again with full computer code available at~\cite{code}. As before, we cannot conclude that the 5 pairs in Table~\ref{table:sharebracket} in the last row are different for large $p$.

\renewcommand{\arraystretch}{1.3}
\begin{table}
\caption{Ineqivalent prime knots that share the KBSM.} \label{table:sharebracket}
\centering
\begin{tabular}{| c | l l l l l |}
\hline
{\bf Space} & \multicolumn{4}{l}{\bf Knots}   & \\
\hline
$L(2,1)$ & \multicolumn{2}{l}{$0_{1}$, $\overline{5_{1}}$, $\overline{5_{5}}$, $\overline{5_{76}}$;} & $\overline{4_{1}}$, $\overline{4_{10}}$; & $\overline{4_{3}}$, $4_{11}$; & $4_{8}$, $\overline{4_{16}}$;\\
& $\overline{4_{23}}$, $5_{1}$; & $\overline{5_{10}}$, $5_{32}$; & $5_{26}$, $5_{27}$; & $\overline{5_{26}}$, $\overline{5_{27}}$; &  $5_{31}$, $5_{60}$; \\
& $5_{36}$, $\overline{5_{60}}$; & $\overline{5_{39}}$, $5_{52}$ & & & \\
 \hline
$L(3,1)$ &  $\overline{1_{1}}$, $\overline{5_{1}}$; & $\overline{2_{1}}$, $\overline{5_{76}}$; & $3_{5}$, $4_{9}$; & $\overline{4_{23}}$, $5_{1}$; & $5_{17}$, $5_{28}$; \\
& $\overline{5_{11}} $, $5_{59}$; & $5_{26}$, $5_{27}$; & $\overline{5_{26}}$, $\overline{5_{27}}$; & $5_{49}$, $5_{56}$ & \\
 \hline
$L(4,1)$ & $2_{2}$, $\overline{5_{1}}$; & $\overline{3_{4}}$, $\overline{5_{76}}$; & $\overline{4_{6}}$, $\overline{5_{77}}$; & $\overline{4_{23}}$, $5_{1}$; & $5_{26}$, $5_{27}$; \\
& $\overline{5_{26}}$, $\overline{5_{27}}$; & $5_{49}$, $\overline{5_{95}}$ & & &\\
 \hline
$L(5,1)$ & $\overline{3_{2}}$, $\overline{5_{1}}$; & $\overline{4_{17}}$, $\overline{5_{76}}$; & $\overline{4_{23}}$, $5_{1}$; & $5_{26}$, $5_{27}$; & $\overline{5_{26}}$, $\overline{5_{27}}$ \\
 \hline
$L(5,2)$ & $3_{6}$, $\overline{5_{1}}$; & $\overline{4_{23}}$, $5_{1}$; & $5_{26}$, $5_{27}$; & $\overline{5_{26}}$, $\overline{5_{27}}$; & $5_{75}$, $\overline{5_{76}}$ \\
 \hline
$L(6,1)$ & $\overline{4_{2}}$, $\overline{5_{1}}$; & $4_{18}$, $5_{3}$; & $\overline{4_{23}}$, $5_{1}$, $5_{72}$;\!\!\! & $5_{26}$, $5_{27}$; & $\overline{5_{26}}$, $\overline{5_{27}}$; \\
& $\overline{5_{68}}$, $\overline{5_{76}}$ & & & & \\
 \hline
$L(p,1), p \leq 12$ &  $4_{23}$, $\overline{5_{1}}$; & $\overline{4_{23}}$, $5_{1}$; & $5_{26}$, $5_{27}$; & $\overline{5_{26}}$, $\overline{5_{27}}$; & $5_{76}$, $\overline{5_{76}}$   \\
 \hline
$L(p,q), p \leq 12, q\geq 2$ &  $4_{23}$, $\overline{5_{1}}$; & $\overline{4_{23}}$, $5_{1}$; & $5_{26}$, $5_{27}$; & $\overline{5_{26}}$, $\overline{5_{27}}$; & $5_{76}$, $\overline{5_{76}}$   \\
 \hline
$L(p,1), p \geq 13$&  $4_{23}$, $\overline{5_{1}}$; & $\overline{4_{23}}$, $5_{1}$; & $5_{26}$, $5_{27}$ \!\!$^\ast$; & $\overline{5_{26}}$, $\overline{5_{27}}$ \!\!$^\ast$; & $5_{76}$, $\overline{5_{76}}$   \\
 \hline
$L(p,q), p \geq 13, q\geq 2$&  $4_{23}$, $\overline{5_{1}}$ \!\!$^\ast$; &\
$\overline{4_{23}}$, $5_{1}$ \!\!$^\ast$; & $5_{26}$, $5_{27}$ \!\!$^\ast$; &\
$\overline{5_{26}}$, $\overline{5_{27}}$ \!\!$^\ast$; &\
$5_{76}$, $\overline{5_{76}}$ \!\!$^\ast$   \\
 \hline
\end{tabular}
\\[0.2cm]
\!\!$^\ast$ Conjecture
\end{table}
\renewcommand{\arraystretch}{1}

\section{Summary} 

For the solid torus $T$ the knot table up to 5 crossings is presented in Appendix A.

For lens spaces $L(p,q)$, $0 < q < p$, $\GCD(p,q)=1$ the knot tables are presented as the knot table for $T$ minus the knots on the left-hand side of the equivalences in Appendix B. 
We have shown that the tables hold for lens spaces $L(p,q)$ with $p \leq 12$, but 2 possible duplicates might appear for $q=1, p \geq 13$ and 5 possible duplicates might appear for $q \geq 2, p \geq 13$.
Such inconclusivenesses are not new to knot theory, for example, the famous Perko pair still causes problems in many classical knot table indexations.

We have found a knot $5_{76}$ for which the $\kbsm$ does not detect its chirality. Furthermore, it is the smallest possible such example.

We have also found a non-invertible link $L(5_{26})$ where methods used by Whitten in \cite{W} fail to detect its non-invertibility. 

The author hopes that the obtained knot tables will serve as a good collection of examples for the ongoing study of knots in $3$-manifolds.

\subsection*{Acknowledgments} The author was supported by the Slovenian Research Agency grant J1-7025. The author would also like to thank M. Mroczkowski and M. Cencelj for their useful insights on the subject.

\section*{Appendix A: The knot table in $\boldsymbol{T}$}\label{app:table}
\begin{center}\begin{tabular}{ccccc}
		\knot{1} & \knot{2} & \knot{4} & \knot{5} & \knot{8}\\
		$0_1$ & $1_{1}$ & $2_{1}$ & ${^\ast}2_{2}$ & $2_{3}$ \\[1em]
		\
		\knot{10} & \knot{12} & \knot{18} & \knot{19} & \knot{21} \\
	    $3_{1}$ &	$3_{2}$ & $3_{3}$ & $3_{4}$ & $3_{5}$ \\[1em]
		\
	   \knot{26} & \knot{28} & \knot{29} & \knot{32} & \knot{34} \\
		 $3_{6}$ & $3_{7}$ & $3_{8}$ & $4_{1}$ & $4_{2}$\\[1em]
		\
		 \knot{40} & \knot{41} &	\knot{43} & \knot{48} & \knot{49} \\
		 $4_{3}$ & $4_{4}$ & $4_{5}$ & $4_{6}$ & $4_{7}$ \\[1em]
		 \
		 \knot{52} & \knot{53} & \knot{64} & \knot{65} &	\knot{67}\\
		  $4_{8}$ & $4_{9}$ & $4_{10}$ & $4_{11}$ & $4_{12}$ \\[1em]
		  \
		  	 \knot{69} & \knot{71} & \knot{73} & \knot{80} &	\knot{82}\\
		 ${^\ast}4_{13}$ & ${^\ast}4_{14}$ & $4_{15}$ & $4_{16}$ &$4_{17}$\\[1em]
		\
	 \knot{84} & \knot{86} & \knot{96} & \knot{97}& 	\knot{100}\\
		 $4_{18}$ & $4_{19}$ & $4_{20}$ & ${^\ast}4_{21}$&$4_{22}$\\[1em]
		 \
		 	 \knot{102} & \knot{110} & \knot{112} & \knot{113} &	\knot{116}\\
		 	 $4_{23}$ & $4_{24}$ & $4_{25}$ & $4_{26}$&${^\ast}4_{27}$\\[1em]
		 	 		\
		 	 		\knot{126} & \knot{128} & \knot{130} & \knot{132}&	\knot{142}\\
		 	 		$5_{1}$ & $5_{2}$ & $5_{3}$ & $5_{4}$&$5_{5}$
	\end{tabular}\end{center}

\begin{center}\begin{tabular}{ccccc}

	 \knot{144} & \knot{150} & \knot{151} & \knot{153}&	\knot{174}\\
		 $5_{6}$ & $5_{7}$ & $5_{8}$ & $5_{9}$&$5_{10}$\\[1em]
		\
	  \knot{176} & \knot{178} & \knot{180} & \knot{182} &	\knot{184}\\
		  $5_{11}$ & $5_{12}$ & $5_{13}$ & $5_{14}$ &	$5_{15}$\\[1em]
		  
		  	  \knot{186} & \knot{188} & \knot{206} & \knot{208} & \knot{210}\\
	  $5_{16}$ & $5_{17}$ & $5_{18}$ & $5_{19}$ & $5_{20}$ \\[1em]
		\
		  \knot{212} & \knot{222} & \knot{223} & \knot{225} & 	\knot{230}\\
		 $5_{21}$ & $5_{22}$ & $5_{23}$ & $5_{24}$ & 	$5_{25}$ \\[1em]
		 		\
	  \knot{231} & \knot{232} & \knot{233} & \knot{254} & 	\knot{255}\\
	 $5_{26}$ & $5_{27}$ & $5_{28}$ & $5_{29}$ & 	$5_{30}$\\[1em]
	 
	 		\
	 		\knot{231} & \knot{232} & \knot{233} & \knot{254} & 	\knot{255}\\
	 		$5_{26}$ & $5_{27}$ & $5_{28}$ & $5_{29}$ & 	$5_{30}$\\[1em]
	 		\
	 		\knot{257} & \knot{262} & \knot{263} & \knot{265}	 & \knot{286}\\
	 		$5_{31}$ & $5_{32}$ & $5_{33}$ & $5_{34}$ & 	$5_{35}$\\[1em]
	 		\
	 		\knot{287} & \knot{289} & \knot{294} & \knot{296} & 	\knot{302}\\
	 		$5_{36}$ & $5_{37}$ & $5_{38}$ & $5_{39}$ & $5_{40}$\\[1em]		
	 		\
	 		\knot{308} & \knot{318} & \knot{320} & \knot{322} & \knot{324}\\
	 		$5_{41}$ & $5_{42}$ & $5_{43}$ & $5_{44}$ & 	$5_{45}$ 
	\end{tabular}\end{center}

\begin{center}\begin{tabular}{ccccc}

		 		  \knot{326} & \knot{328} & \knot{330} & \knot{332} & \knot{350} \\
	 $5_{46}$ & $5_{47}$ & $5_{48}$ & $5_{49}$ & $5_{50}$ \\[1em]
		\
 \knot{351} & \knot{354} & \knot{358} & \knot{362} & 	\knot{363}\\
	 $5_{51}$ & $5_{52}$ & $5_{53}$ & $5_{54}$ & 	$5_{55}$\\[1em]
		\
	 \knot{385} & \knot{390} & \knot{391} & \knot{393} & 	\knot{395}\\
	 $5_{56}$ & $5_{57}$ & $5_{58}$ & $5_{59}$ & $5_{60}$\\[1em]
\
			 \knot{397} & \knot{398} & \knot{400} & \knot{404} & 	\knot{406}\\
		 $5_{61}$ & $5_{62}$ & $5_{63}$ & $5_{64}$ & 	$5_{65}$ \\[1em]
		\
	 \knot{411} & \knot{424} & \knot{426} & \knot{428} & \knot{430}\\
 $5_{66}$ & $5_{67}$ & $5_{68}$ & $5_{69}$ & $5_{70}$\\[1em]
		\
		 \knot{446} & \knot{448} & \knot{450} & \knot{456} &\\
		 $5_{71}$ & $5_{72}$ & $5_{73}$ & $5_{74}$ &
	\end{tabular}\end{center}
\vspace{1em}
$\ast$ -- amphichiral knot


\section*{Appendix B: Equivalences of knots in lens spaces}\label{app:eq}
\renewcommand\arraystretch{1.2}
\begin{table}[H]
\centering

 \small

\caption{Equivalences 	of knots in lens spaces.}
\begin{tabular}{ |c | r@{$\;\;\;\sim\;\;\;$}l !{\color{blu}\vrule} r@{$\;\;\;\sim\;\;\;$}l !{\color{blu}\vrule} r@{$\;\;\;\sim\;\;\;$}l |} \hline
{\bf Space} & \multicolumn{6}{c|}{ {\bf Equivalences}} \\ \hline
$\mathbf{L(2,1)}$ & $1_{1}$, $5_{74}$ & $O$ &  $4_{16}$ & $\overline{3_{4}}$ & $\overline{5_{2}}$ & $\overline{4_{6}}$\\
& $5_{24}$ & $\clubsuit$ & $2_{3}$, $4_{22}$, $4_{23}$, $5_{76}$, $5_{78}$ & $0_1$& $\overline{5_{41}}$, $\overline{5_{90}}$ & $\overline{3_{3}}$ \\
& $4_{27}$ & $\overline{4_{7}}$ & $5_{54}$ & $5_{32}$ & $3_{6}$, $\overline{3_{8}}$ & $\overline{1_{1}}$ \\
& $5_{85}$ & $\overline{3_{6}}$ & $5_{87}$ & $4_{10}$ & $\overline{5_{51}}$, $\overline{5_{92}}$ & $\overline{5_{32}}$ \\
& $3_{7}$ & $2_{1}$ & $4_{1}$, $5_{88}$ & $3_{8}$ & $5_{39}$ & $\overline{4_{11}}$ \\
& $\overline{5_{87}}$ & $\overline{5_{43}}$ & $4_{26}$, $5_{75}$, $5_{81}$ & $\overline{2_{1}}$ & $\overline{4_{22}}$ & $\overline{4_{2}}$ \\
& $\overline{4_{25}}$ & $\overline{4_{17}}$ & $\overline{5_{78}}$ & $\overline{5_{68}}$ & $4_{24}$ & $\overline{2_{3}}$ \\
& $5_{93}$ & $4_{6}$ & $5_{89}$ & $4_{20}$ & $5_{77}$ & $0_{1}\#\clubsuit$ \\
& $\overline{3_{7}}$, $\overline{5_{94}}$ & $\overline{3_{2}}$ & $5_{42}$ & $4_{7}$ & $5_{5}$ & $\overline{4_{23}}$ \\
& $5_{91}$ & $5_{79}$ & $4_{25}$, $5_{71}$ & $3_{3}$ & $\overline{5_{4}}$, $\overline{5_{91}}$ & $4_{9}$ \\
& $\overline{5_{38}}$, $5_{67}$ & $\overline{4_{26}}$ & $\overline{5_{93}}$ & $5_{80}$ & $\overline{5_{95}}$ & $3_{5}$ \\
& $5_{3}$, $5_{23}$ & $\overline{4_{8}}$ & $\overline{5_{74}}$ & $\overline{5_{6}}$ & \multicolumn{2}{c |}{}\\
 \hline \hline 
$\mathbf{L(3,1)}$ &   $2_{3}$, $5_{71}$ & $O$ &$1_{1}$, $3_{6}$, $5_{74}$, $5_{75}$ & $0_1$ & $\overline{5_{2}}$, $5_{77}$ & $\overline{4_{7}}$ \\
& $\overline{5_{4}}$ & $5_{3}$ & $4_{23}$, $4_{24}$, $5_{24}$ & $\overline{1_{1}}$ & $4_{9}$ & $3_{5}$ \\
& $5_{29}$ & $4_{15}$ & $5_{33}$, $5_{82}$ & $5_{22}$ & $4_{22}$ & $2_{1}$ \\
& $4_{26}$ & $\overline{3_{4}}$ & $\overline{5_{90}}$ & $\overline{4_{11}}$ & $\overline{5_{41}}$ & $5_{23}$ \\
& $\overline{3_{8}}$ & $2_{2}$ & $5_{5}$ & $3_{8}$ & $4_{25}$ & $4_{16}$ \\
& $5_{55}$ & $5_{34}$ & $5_{76}$ & $\overline{2_{1}}$ & $\overline{4_{2}}$ & $\overline{3_{7}}$ \\
& $5_{67}$ & $\overline{4_{17}}$ & $\overline{5_{91}}$ & $5_{39}$ & $5_{78}$, $5_{81}$ & $3_{3}$ \\
& $4_{20}$ & $4_{6}$ & $\overline{5_{6}}$ & $\overline{4_{22}}$ & $5_{93}$ & $5_{40}$ \\
& $5_{85}$ & $\overline{2_{3}}$ & $5_{88}$ & $4_{7}$ & $\overline{5_{68}}$ & $\overline{4_{25}}$ \\
& $5_{87}$ & $5_{42}$ & $3_{7}$ & $3_{1}$ & $\overline{5_{94}}$ & $4_{8}$ \\
& $\overline{5_{38}}$ & $4_{27}$ &\multicolumn{2}{c !{\color{blu}\vrule}}{} & \multicolumn{2}{c |}{} \\
\hline \hline

$\mathbf{L(4,1)}$ &  $4_{24}$, $2_{3}$ & $0_1$  & $\overline{3_{8}}$ & $\overline{3_{2}}$ & $5_{78}$ & $4_{16}$ \\
& $5_{77}$ & $\overline{5_{38}}$ & $5_{75}$, $5_{85}$ & $\overline{1_{1}}$ & $5_{76}$ & $\overline{3_{4}}$ \\
& $5_{39}$ & $4_{18}$ & $5_{72}$ & $5_{40}$ & $5_{74}$ & $2_{1}$ \\
& $3_{6}$ & $O$ & $4_{26}$ & $\overline{4_{17}}$ & $\overline{5_{90}}$ & $\overline{5_{43}}$ \\
& $4_{23}$ & $2_{2}$ & $4_{1}$ & $3_{7}$ & $5_{67}$ & $4_{25}$ \\
& $5_{88}$ & $5_{79}$ & $4_{22}$ & $3_{1}$ & $\overline{5_{6}}$ & $\overline{3_{7}}$ \\
& $5_{71}$ & $5_{24}$ & $\overline{5_{94}}$ & $5_{80}$ &\multicolumn{2}{c |}{} \\
 \hline
 
 $\mathbf{L(5,1)}$ & $5_{85}$, $3_{6}$ & $0_1$ & $5_{74}$ & $3_{1}$ & $\overline{4_{2}}$ & $\overline{3_{8}}$ \\
 & $4_{24}$ & $O$ & $2_{3}$ & $1_{1}$ & $4_{23}$ & $\overline{3_{2}}$ \\
 & $4_{22}$ & $4_{1}$ & $\overline{5_{68}}$ & $4_{26}$ & $5_{75}$ & $2_{2}$ \\
 & $5_{5}$ & $3_{7}$ & $5_{76}$ & $\overline{4_{17}}$ & $5_{78}$ & $5_{67}$ \\
 \hline \hline 
 $\mathbf{ L(5,2)}$ & $4_{24}$ & $0_1$ & $2_{3}$ & $1_{1}$ & $4_{23}$ & $3_{6}$ \\
 & $5_{76}$ & $5_{75}$  &\multicolumn{2}{c !{\color{blu}\vrule}}{}&\multicolumn{2}{c |}{}\\
 \hline \hline

 $\mathbf{L(6,1)}$ & $4_{24}$ & $0_1$ & $3_{6}$ & $1_{1}$ & $\overline{5_{6}}$ & $\overline{3_{8}}$ \\
 & $4_{23}$ & $\overline{4_{2}}$ & $5_{85}$ & $O$ & $5_{75}$ & $\overline{3_{2}}$ \\
 & $5_{74}$ & $4_{1}$ & $5_{5}$ & $4_{22}$ & $5_{76}$ & $\overline{5_{68}}$ \\
 \hline \hline

 $\mathbf{L(7,1)}$ & $5_{85}$ & $0_1$ & $3_{6}$ & $2_{3}$ & $\overline{5_{6}}$ & $4_{23}$  \\
 & $4_{24}$ & $1_{1}$ & $5_{75}$ & $\overline{4_{2}}$ & $5_{74}$ & $5_{5}$ \\
 \hline
\end{tabular}
\end{table}
\clearpage

\begin{table}[H]
\centering
\begin{tabular}{ |c | r@{$\;\;\;\sim\;\;\;$}l !{\color{blu}\vrule} r@{$\;\;\;\sim\;\;\;$}l !{\color{blu}\vrule} r@{$\;\;\;\sim\;\;\;$}l |} \hline
{\bf Space} & \multicolumn{6}{c|}{ {\bf Equivalences}} \\ \hline

 $\mathbf{L(7,2)}$  & $3_{6}$ & $2_{3}$ & $5_{75}$ & $4_{24}$  &\multicolumn{2}{c |}{} \\
 \hline \hline
$\mathbf{L(8,1)}$ & $5_{85}$ & $1_{1}$ & $4_{24}$ & $2_{3}$ & $5_{75}$ & $\overline{5_{6}}$ \\
 \hline \hline
 $\mathbf{L(8,3)}$ &\multicolumn{2}{c !{\color{blu}\vrule}}{/} & \multicolumn{2}{c !{\color{blu}\vrule}}{} & \multicolumn{2}{c |}{} \\
  \hline \hline
$\mathbf{L(9,1)}$  & $5_{85}$ & $2_{3}$ & $4_{24}$ & $3_{6}$ &\multicolumn{2}{c |}{}  \\
  \hline \hline
$\mathbf{L(9,2)}$ & $4_{24}$ & $3_{6}$   &  \multicolumn{2}{c !{\color{blu}\vrule}}{}  & \multicolumn{2}{c |}{}  \\

$\mathbf{L(10,1)}$  & $5_{85}$ & $3_{6}$  & \multicolumn{2}{c !{\color{blu}\vrule}}{\;\;\;\;\;\;\;\;\;\;\;\;\;\;\;\;}  & \multicolumn{2}{c |}{\;\;\;\;\;\;\;\;\;\;\;\;\;\;\;\;}  \\
  \hline \hline
$\mathbf{L(10,3)}$ & \multicolumn{2}{c !{\color{blu}\vrule}}{/} & \multicolumn{2}{c !{\color{blu}\vrule}}{} & \multicolumn{2}{c |}{} \\
  \hline \hline
$\mathbf{ L(11,1)}$ & $5_{85}$ & $4_{24}$ &  \multicolumn{2}{c !{\color{blu}\vrule}}{}  & \multicolumn{2}{c |}{}  \\
  \hline \hline
$\mathbf{L(11,2)}$ & $5_{85}$ & $4_{24}$  & \multicolumn{2}{c !{\color{blu}\vrule}}{}  & \multicolumn{2}{c |}{}  \\
 \hline \hline
$\mathbf{L(11,3)}$ &  \multicolumn{2}{c !{\color{blu}\vrule}}{/} & \multicolumn{2}{c !{\color{blu}\vrule}}{} & \multicolumn{2}{c |}{} \\
 \hline \hline
 $\mathbf{L(12,1)}$ &  \multicolumn{2}{c !{\color{blu}\vrule}}{/} & \multicolumn{2}{c !{\color{blu}\vrule}}{} & \multicolumn{2}{c |}{} \\
 \hline \hline
  $\mathbf{L(12,5)}$ & \multicolumn{2}{c !{\color{blu}\vrule}}{/} & \multicolumn{2}{c !{\color{blu}\vrule}}{} & \multicolumn{2}{c |}{} \\
 \hline \hline
  $\mathbf{L(p,q), p\geq 13}$  & \multicolumn{2}{c !{\color{blu}\vrule}}{/$^\ast$} & \multicolumn{2}{c !{\color{blu}\vrule}}{} & \multicolumn{2}{c |}{} \\
 \hline 

 \end{tabular}
\end{table}
\renewcommand\arraystretch{1.0}

\begin{tabular}{l l@{ -- }l}
Legend: & $O$ & affine unknot \\
& $\clubsuit$ & affine trefoil \\
& $^\ast$ & conjecture
\end{tabular}





\end{document}